\documentclass[11pt]{article}
\usepackage{amsmath,amsfonts,amssymb}

\newcommand{\ux}{\underline x}
\newcommand{\uy}{\underline y}
\newcommand{\Abar}{\widetilde A}
\newcommand{\Span}[1]{\left<#1\right>}
\newcommand{\dd}{\mathrm d}
\newcommand{\tM}{\,{}^t\!M} 
\newcommand{\tP}{\,{}^t\!P} 

\newcommand{\half}{\frac12}
\newcommand{\fie}{\varphi}

\newcommand{\om}{\omega}

\newcommand{\Ga}{\Gamma}
\newcommand{\C}{\mathbb C}
\newcommand{\Q}{\mathbb Q}

\newcommand{\PP}{\mathbb P}
\newcommand{\Oh}{\mathcal O}
\newcommand{\sI}{\mathcal I}

\newcommand{\iso}{\cong}
\newcommand{\ot}{\leftarrow}
\newcommand{\into}{\hookrightarrow}
\newcommand{\onto}{\twoheadrightarrow}
\DeclareMathOperator{\codim}{codim}

\DeclareMathOperator{\rank}{rank}
\DeclareMathOperator{\Basis}{Basis}

\DeclareMathOperator{\sExt}{\mathcal{E}\mathit{xt}}

\DeclareMathOperator{\sHom}{\mathcal{H}\mathit{om}}
\DeclareMathOperator{\Pf}{Pf}

\DeclareMathOperator{\Basket}{Basket}
\DeclareMathOperator{\Numerator}{Numerator}

\newtheorem{thm}{Theorem}[section]

 \newtheorem{exa}[thm]{Example}
 
 \newtheorem{rmk}[thm]{Remark}

\numberwithin{equation}{section}

\title{Examples of Type IV unprojection}
\author{Miles Reid}
\date{}
\begin{document}
\maketitle

 \begin{abstract}
I show that $\PP(2,3)$ has an embedding $\PP(2,3)\iso\Ga\subset
\PP(4,5,6,9)$ whose image $\Ga$ is contained in a quasismooth K3
hypersurface $X_{24}\subset\PP(4,5,6,9)$. The pair $\Ga\subset X_{24}$
unprojects to the codimension~4 K3 surface $Y\subset\PP(4,5,5,6,7,8,9)$
with
 \begin{align*}
 \Basket &= [\,\textstyle{\half(1,1), \frac15(1,4), \frac15(2,3), \frac19(4,5)
 }\,] .\\
 \Numerator &= 1-t^{12}-t^{13}-2t^{14}-2t^{15}-2t^{16}-t^{17} \\
 & \kern2cm +t^{19}+2t^{20}+3t^{21}+4t^2+3t^{23}+\cdots
 \end{align*}
(Alt{\i}nok4(111) in the Magma K3 database). The local coordinates at the
third centre $P_3=\frac15(2,3)$ of $Y$ are of weight 7 and 8 (rather than
2 and 3), so both are eliminated by the projection from $P_3$. Together with
other examples, this gives substance to Type~IV unprojections. Several more
cases of Type~IV unprojections are known up to codimension~5 or~6. The paper
also contains some Magma programming routines suitable as exercises for
babies.
 \end{abstract}

 \section{Introduction}
 Unprojection makes new Gorenstein rings out of old. The idea is to use
the adjunction formula in Serre--Grothendieck duality to construct
rational functions with poles along a given divisor. Adjoining these functions
is an analog of Castelnuovo's contractibility criterion in minimal models
of surfaces. Type~I unprojections were introduced by Kustin and Miller
\cite{KM} and later by Papadakis and Reid \cite{PR}. They also figure
prominently in birational geometry, see Corti, Pukhlikov and Reid
\cite{CPR}, and \cite{Ki}, and in that context can be traced back to work
on the Cremona group in the late 19th century (Geiser involutions, see
\cite{CPR}, 2.6.3.). Unprojections become progressively more complicated as the
divisor to be unprojected gets further from projectively Gorenstein, so
that several unprojection variables have to be adjoined at the same time.
Type~II projections first occurred in \cite{CPR}, 4.10--4.12 and 7.3 (in
which context they are related to Bertini involutions), and are studied
in more detail in \cite{Ki}, Section~9. The typical case of Type~III
projections is the classical projection of a Fano 3-fold
$V_{2g-2}\subset\PP^{g+1}$ from a line (compare \cite{Ki}, Section~9 and
Example~9.16).

 This paper introduces Type~IV unprojections. I don't intend to give a
formal definition here. The defining feature of Type~IV is that the
coordinate ring of the divisor to be unprojected is a ring $A=\C[\Ga]$
whose normalisation $\Abar$ is Gorenstein, and such that the quotient
$\Abar/A$ is an Artinian module needing 2 linearly independent generators
over $A$. An ideal case would be
 \begin{equation}
 A=\C[\Ga]:=\C[u^2,uv,v^2,u^3,u^2v,uv^2,v^3],
 \quad\hbox{with}\quad \Abar=\C[u,v].
 \end{equation}
Here $\Abar$ needs $1,u,v$ as generators over $A$, and $u,v$ are
linearly independent. I do not treat this case here, because any example
with it would have codimension $\ge4$ unprojecting to codimension
$\ge7$; see however Section~\ref{sec!det}

The Type~IV unprojections studied in this paper add 3 new variables as
generators to the ambient ring, so potentially increase the minimal number
of generators by~3 (hence also the codimension). But this increase is
masked if some of the old variables become expressible in terms of the new;
this frequently happens in small codimension. The K3 database in Magma
(export~2.8, \cite{Ma}) reports 5 families of K3 surfaces in codimension~4
having a Type~IV unprojection to codimension~1 (see \ref{ssec!K3DB} and
compare \cite{ABR}); there are also 6 families in codimension~5 having a
Type~IV unprojection to codimension~2 and several similar things in
codimension~6.

 \section{Embedding of $\Ga$ and quasismoothness of $X_{24}$} \label{sec!emb}
Write $u,v$ for coordinates on $\PP(2,3)$ and $x,y,z,t$ for $\PP(4,5,6,9)$.
Up to change of coordinates, the general map
$i\colon\PP(2,3)\to\PP(4,5,6,9)$ is
 \begin{equation}
 i:(u,v)\mapsto(u^2,uv,u^3+v^2,u^3v+v^3).
 \label{eq!i}
 \end{equation}
Choosing the 4th entry $t=u^3v+v^3$ (so that $t=vz$) tidies up later
calculations. It is an embedding because monomials in
$(u^2,uv,u^3+v^2,u^3v+v^3)$ span the vector space $H^0(\PP(2,3),\Oh(i))$
for all $i\gg0$ (in fact for every $i\ge13$, see \ref{ssec!smooth}).

The image $\Ga=i(\PP(2,3))$ is defined by
 \begin{equation}
 \begin{array}{rl}
 f_{18}:& x^3z+xyt-z^3+t^2, \\
 g_{19}:& x^2yz-xzt+y^2t, \\
 h_{20}:& x^5-x^2z^2+2xy^2z-y^4,
 \end{array}
 \qquad
 \begin{array}{rl}
 q_{21}:& x^3t-xyz^2+y^3z, \\
 q_{22}:& -xt^2+y^2z^2, \\
 q_{23}:& -x^2zt+yz^3-yt^2
 \end{array}
 \label{eq!Ga}
 \end{equation}
(see \ref{ssec!smooth}). For reasons of degree, only the first 3 of
these can take part in the equation of $X_{24}$: set
 \begin{equation}
 F_{24}=xh_{20}+yg_{19}+zf_{18}.
 \label{eq!F}
 \end{equation}
One checks by brute force (see \ref{ssec!smooth}) that $F$ defines a
quasismooth hyper\-surface $X_{24}\subset\PP(4,5,6,9)$.

 \section{Resolution of $\C[u,v]$ as a $\C[x,y,z,t]$ module}
 Write $A=\C[x,y,z,t]$ for the homogeneous coordinate ring of $\PP(4,5,6,9)$
and $\C[\Ga]=A/I_\Ga$ for that of $\Ga\subset\PP(4,5,6,9)$; by construction,
its normalisation is the coordinate ring $\C[u,v]$ of $\PP(2,3)$. Thus the
embedding $i$ of (\ref{eq!i}) makes $\C[u,v]$ a module over the polynomial ring
$\C[x,y,z,t]$, with
 \begin{equation}
 \begin{gathered}
 x\cdot1=u^2,\quad y\cdot1=uv,\quad z\cdot1-x\cdot u=v^2,
 \quad x\cdot u=u^3, \\
 x\cdot v=u^2v, \quad y\cdot v=uv^2, \quad z\cdot v-x^2\cdot1=v^3,
 \quad \hbox{etc.}
 \end{gathered}
 \label{eq!elim}
 \end{equation}
This already shows that $\C[u,v]$ is generated as a $\C[x,y,z,t]$-module
by the 3 elements $(1,u,v)$. The relations between these generators are
 \begin{equation}
 (v,u,1)M=0,
 \end{equation}
where
 \begin{equation}
 M= \begin{pmatrix}
 -x & -y & -z & -t & 0 & 0 \\
 y & z & 0 & -xz & t & x^2 \\
 0 & -x^2 & t & z^2 & -yz & -xz+y^2
 \end{pmatrix}
 \label{eq!M}
 \end{equation}
(see \ref{ssec!res}). The $3\times6$ matrix $M$ is homogeneous with entries
of weight
 \begin{equation}
 \begin{pmatrix}
 4 & 5 & 6 & 9 & 8 & 7 \\
 5 & 6 & 7 & 10 & 9 & 8 \\
 7 & 8 & 9 & 12 & 11 & 10
 \end{pmatrix}
 \end{equation}
and satisfies
 \begin{equation}
 MJ\tM = 0, \quad\hbox{where}\quad
 J=\begin{pmatrix} 0 & I_3 \\ -I_3 & 0 \end{pmatrix}.
 \end{equation}

Now $\C[u,v]$ is a Gorenstein module of codim~2 over $\C[x,y,z,t]$,
and one checks that its resolution is
 \begin{equation}
 \C[u,v] \longleftarrow L_0 \xleftarrow{\,M\,} L_1
 \xleftarrow{\,J\tM\,} L_2 \ot 0,
 \label{eq!Gor}
 \end{equation}
where
 \begin{align*}
 L_0 &= A\oplus A(-2)\oplus A(-3), \\[4pt]
 L_1 &= A(-7)\oplus A(-8)\oplus A(-9)\oplus A(-12)\oplus A(-11)\oplus A(-10),
 \\[4pt]
 L_2 &= A(-16)\oplus A(-17)\oplus A(-19).
 \end{align*}

 \section{Maps between complexes}
 I find the unprojection variables by comparing the resolutions of
$\C[X]=A/(F_{24})$ and $\C[u,v]$. For this, consider first the diagram
including part of the resolution of $\C[\Ga]$:
 \begin{equation}
 \renewcommand{\arraystretch}{1.5}
 \begin{matrix}
 \C[X] & \ot & A & \ot & A(-24) & \ot & 0 \\
 \big\downarrow && \Vert && \big\downarrow \\
 \C[\Ga] & \ot & A & \ot & K_1 & \ot & \cdots \\
 \bigcap && \bigcap && \big\downarrow \\
 \C[u,v] & \ot & L_0 & \ot & L_1 & \ot & L_2 & \ot & 0
 \end{matrix}
 \label{eq!comps}
 \end{equation}
Here $K_1=A(-18)\oplus\cdots\oplus(-23)$ is the free module corresponding
to the 6 generators (\ref{eq!Ga}) of the ideal $I_\Ga$. The first downarrow
$A(-24)\to K_1$ is the matrix $(z,y,x,0,0,0)$ that expresses $F_{24}$ as
the combination (\ref{eq!F}) of $f_{18},g_{19},h_{20}$. For any $f\in I_\Ga$,
the column vector $(0,0,f)\in\ker\{L_0\to\C[u,v]\}$, and is therefore hit by
some element of $L_1$. The second downarrow is the $6\times6$ matrix
expressing the 6 generators (\ref{eq!Ga}) of $K_1$ as linear combinations
of the 6 columns of $M$. This matrix is presumably implicit in the Gr\"obner
basis calculation that eliminates $u,v$. Because $F_{24}$ is a linear
combination of only the first three generators (\ref{eq!F}), I only need
the corresponding three columns:
 \begin{equation}
 \begin{pmatrix}
 0&0&0 \\
 0&0&0 \\
 f_{18}&g_{19}&h_{20} \\
 \end{pmatrix} = MN,
 \label{eq!MN}
 \end{equation}
where
 \begin{equation}
 N = \begin{pmatrix}
 0 & 0 & x^2y+xt \\
 -xz & 0 & -x^3 \\
 t+xy & 0 & 0 \\
 -z & 0 & -x^2 \\
 0 & -x^2 & -xy \\
 0 & t & -y^2
 \end{pmatrix}.
 \end{equation}
Composing gives
 \begin{equation}
 \renewcommand{\arraystretch}{1.5}
 \begin{matrix}
 \C[X] & \ot & A & \ot & A(-24) & \ot & 0 \\
 \big\downarrow && \bigcap && \big\downarrow \\
 \C[u,v] & \ot & L_0 & \ot & L_1 & \ot & L_2 & \ot & 0,
 \end{matrix}
 \label{eq!cxs}
 \end{equation}
where the downarrow takes the basis $F_{24}$ to
 \begin{equation}
 N \begin{pmatrix}
 z \\ y \\ x
 \end{pmatrix} =
 \begin{pmatrix}
 x^3y+x^2t \\
 -x^4-xz^2 \\
 xyz+zt \\
 -x^3-z^2 \\
 -2x^2y \\
 -xy^2+yt
 \end{pmatrix}
 \in L_1(24).
 \label{eq!N24}
 \end{equation}

This allows me to start writing down the unprojection $Y$ of $\Ga$ in $X$.
The theory is similar to that of Kustin and Miller \cite{KM}, and Papadakis
and Reid \cite{PR} and \cite{Ki}: the adjunction formula
$\om_\Ga=\sExt^1_X(\Oh_\Ga,\om_X)$ gives an exact sequence
 \begin{equation}
 0\to\om_X \to \sHom_{\Oh_X}(\sI_\Ga,\om_X)\to \om_\Ga\to 0.
 \end{equation}
Note that I abuse notation by writing $X,\Ga$ in place of the punctured
affine cones over them. The modules and Homs between them are really the
Serre modules of sheaves on $X,\Ga$. For example, $\Ga\iso\PP^1$, but
as $\PP(2,3)$ it has $\om_\Ga=\Oh(-5)$, since on the punctured affine cone
over $\Ga$, the dualising sheaf twisted by $\Oh(5)$ is based by
$\dd u\wedge\dd v$.

Moreover since $\C[\Ga]\into\C[u,v]$ is an isomorphism outside the origin,
the dualising module $\om_{\C[\Ga]}=\om_{\C[u,v]}$ is isomorphic to
$\C[u,v]$, so is free of rank~1 as a $\C[u,v]$ module; however, as a
module over $\C[\Ga]$ or $\C[X]$ or $\C[x,y,z,t]$, it needs 3 generators,
corresponding to $1,u,v$. I used the fact that $\C[u,v]$ is a Gorenstein
module to derive the resolution (\ref{eq!Gor}); it is best not to enquire
too closely about the resolution of $\C[\Ga]$.

Now the homogeneous coordinate ring $\C[Y]$ of the unprojected variety $Y$
is obtained from $\C[X]$ by adjoining rational functions $s_0,s_1,s_2$ with
poles along $\Ga$, viewed as homomorphisms $\sI_\Ga\to\om_X\iso\Oh_X$.
Here $s_0\in\sHom(\sI_\Ga,\om_X(5))$ is the rational form with pole along
$\Ga$ whose Poincar\'e residue is a basis of $\om_\Ga(5)\iso\Oh_{\PP(2,3)}$.
As in the proof of Castelnuovo's criterion, the graded ring with $s_0$
adjoined already contracts $\Ga$ to a point. The elements $s_1,s_2$ that map
to $u,v$ times the generator of $\om_\Ga(5)$ under
$\sHom_{\Oh_X}(\sI_\Ga,\om_X)\onto\om_\Ga$ are required to make the image
projectively normal.

I choose $s_0$ to map to $1\in\C[u,v]$, and $s_1,s_2$ to $u\cdot1,v\cdot1$
respectively; since $\om_{\PP(2,3)}=\Oh(-5)$ and $\om_X=\Oh_X$,
 \begin{equation}
 \deg s_0=5,\quad\hbox{and}\quad \deg s_1,s_2=7,8.
 \end{equation}
The {\em linear} relations between $s_0,s_1,s_2$ come from (\ref{eq!cxs}).
Namely (essentially as in \cite{KM}, see also \cite{PR}, end of Section~1),
$s_0,s_1,s_2$ correspond to the generators of $L_2$ or the rows of $M$, and
$J\tM\colon L_2\to L_1$ must map them to the image of $F_{24}$ under
the downarrow. This gives the equations
 \[
 J\tM  \begin{pmatrix} s_2 \\ s_1 \\ s_0 \end{pmatrix} =
 N \begin{pmatrix}
 z \\ y \\ x
 \end{pmatrix},
 \]
or spelled out:
 \begin{equation}
 \begin{pmatrix}
 -t & - xz & z^2 \\
 0 & t & - yz \\
 0 & x^2 & - xz+y^2 \\
 x & - y & 0 \\
 y & - z & x^2 \\
 z & 0 & - t 
 \end{pmatrix}
 \begin{pmatrix} s_2 \\ s_1 \\ s_0 \end{pmatrix} =
 \begin{pmatrix}
x^3y+x^2t \\
-x^4-xz^2 \\
xyz+zt \\
-x^3-z^2 \\
-2x^2y \\
-xy^2+yt
 \end{pmatrix}.
 \label{eq!sp}
 \end{equation}
I write $R_1,\dots,R_6$ for these relations. The equations $f,g,h$ of
$\Ga$ appear naturally $3\times3$ minors of $J\tM$ (or $2\times2$ minors
of defective $2\times3$ blocks with zeros down one column) so that
$R_1,\dots,R_6$ can be solved by linear algebra to give expressions for
$(f,g,h)$ times $s_2,s_1,s_0$. For example
 \begin{align*}
fs_0&=-3x^3yz+x^2y^3-x^2zt-yt^2,\\
gs_0&=x^6+x^3z^2+xyzt+zt^2,\\
hs_0&=-x^5y+2x^2yz^2-xy^3z+xz^2t-y^2zt.\\
 \end{align*}
This makes explicit that $s_0$ is a homomorphism $\sI_\Ga\to\Oh_X$. But
it is much more useful to use $R_1,\dots,R_6$ as they are, without taking
determinants.

 \section{The quadratic relations for $s_1^2,s_1s_2,s_2^2$} \label{sec!quad}
 Quadratic relations in $s_1,s_2$ must exist for several reasons. One is
that the ring extension $\C[X]\subset\C[Y]$ can be obtained by adjoining
$s_0$ first, followed by the normalisation $\C[X][s_0]\subset\C[Y]$; the
second step is integral, and one guesses that $\C[Y]$ is generated by
$1,s_1,s_2$ as a module over $\C[X][s_0]$. Or, since $s_1,s_2$ correspond
in some sense to $us_0$ and $vs_1$ on taking Poincar\'e residue to $\Ga$,
and since, restricted to $\Ga$, $x=u^2,y=uv,z=v^2+xu$, there must be
relations $S_0,S_1,S_2$ saying that
 \begin{equation}
 \renewcommand{\arraystretch}{1.3}
 \begin{array}{rrl}
 s_1^2-xs_0^2 &\in&\Span{xys_0,ts_0,zs_2,x^2z,xy^2,yt}, \\
 s_1s_2-ys_0^2 &\in&\Span{xzs_0,y^2s_0,x^2s_1,xyz,y^3,zt}, \\
 s_2^2-zs_0^2+xs_0s_1 &\in&\Span{yzs_0,xys_1,ts_1, x^4,xz^2,y^2z}, \\
 \end{array}
 \end{equation}
where the right-hand side just lists alll monomials of degree 14, 15 and 16.

To find these relations explicitly, I proceed as follows. Since $MJ\tM=0$
and the first two rows of $MN$ are zero by (\ref{eq!MN}), the first two rows
of $M$ provide two syzygies between the 6 rows $R_1,\dots,R_6$ of (\ref{eq!sp}).
I massage the second very slightly to make them both 4-term syzygies:
 \begin{equation}
 \begin{array}{rcll}
 xR_1+yR_2+zR_3+tR_4 &\equiv&0,&\hbox{and}\\ [4pt]
 yR_1+z(R_2-xR_4)+tR_5+x^2R_6 &\equiv&0.
 \end{array}
 \end{equation}
This suggests realising the relations as $4\times4$ Pfaffians of
$5\times5$ skew matrixess, as follows:
 \begin{equation}
 \renewcommand{\arraycolsep}{0.33em}
 \begin{pmatrix}
 x&y&z&t \\
 &z&s_1&ys_0-xz \\
 && s_2+x^2 & -xs_1+zs_0 \\
 &&& x^3
 \end{pmatrix},
\quad
 \begin{pmatrix}
 y&z&x^2&t \\
 & s_0+y&s_1&s_2 \\
 &&s_2+x^2&-xy \\
 &&& xs_1-zs_0
 \end{pmatrix}.
 \end{equation}
My convention is to write only the upper triangular terms
$a_{12},a_{13},\dots,\linebreak[2] a_{45}$ of a skew matrix. The Pfaffians
are defined by $\Pf_{ij.kl}=a_{ij}a_{kl}-a_{ik}a_{jl}+a_{il}a_{jk}$ (up to
$\pm$).

The second matrix has first four Pfaffians $R_5,R_6,R_2-xR_4,R_1$, and the
fifth is one of the required quadratic relations:
 \begin{equation}
 S_2:\ s_2^2-zs_0^2+xs_0s_1=yzs_0-2xys_1-x^2s_2.
 \end{equation}
The first matrix gives $R_1,R_2,R_3,R_4$ and
 \begin{equation}
 \Pf_{23.45}= xs_1^2-zs_0s_1+ys_0s_2-x^2ys_0-xzs_2,
 \label{eq!xs1^2}
 \end{equation}
which is not quite what we want. However, subtracting $s_0R_5$ gives
a relation that is identically divisible by $x$, and taking out the
factor gives
 \begin{equation}
 S_0:\ s_1^2-xs_0^2=xys_0+zs_2.
 \label{eq!s1^2}
 \end{equation}
 
The equation for $s_1s_2$ can be found by constructing an explicit syzygy,
for example by taking $s_1R_4+yS_0+zR_5$, that has $xs_1s_2$ as one term,
and all terms of which are divisible by $x$. However, this is again a
4-term syzygy, once again best understood by realising the relations as
Pfaffians:
 \begin{equation}
 \begin{pmatrix}
 x&y&z&s_1 \\
 &z&s_1&s_2 \\
 && s_2+x^2 & -xs_0-2xy \\
 &&& -s_0^2-ys_0
 \end{pmatrix}.
 \end{equation}
Thus
 \begin{equation}
 S_1:\ s_1s_2-ys_0^2=(y^2-xz)s_0-x^2s_1-2xyz.
 \label{eq!s1s2}
 \end{equation}
In fact, all the relations $R_1,\dots,R_5,S_0,S_1,S_2$ (but, annoyingly,
not $R_6$) can be written together as the $4\times4$ Pfaffians of the
following $6\times6$ matrix:
 \begin{equation}
 \begin{pmatrix}
 x&y&z&s_1&t \\
 &z&s_1&s_2&ys_0-xz \\
 && s_2+x^2&-xs_0-2xy&-xs_1+zs_0 \\
 &&& -s_0^2-ys_0&x^3 \\
 &&&&s_0s_2+x^2s_0+x^2y+xt
 \end{pmatrix}.
 \end{equation}

 \section{Another example} \label{sec!det}
 Consider the following two families of K3 surfaces:
 \begin{enumerate}
 \renewcommand{\labelenumi}{(\alph{enumi})}
 \item The codimension 2 complete intersection
$X_{6,6}\subset\PP(2,2,2,3,3)$, known to its friends as Fletcher2(14),
with $9\times\half(1,1)$ singularities on the plane $\PP(2,2,2)=\PP^2$.
 \item The codimension~5 K3 surface $Y\subset\PP(2,2,2,2,3,3,3,3)$ in
symmetric determinantal format, with 10 nodes. Write
$x_1,\dots,x_4,y_1,\dots,y_4$ for coordinates, and let $M$ be a general
$4\times4$ symmetric matrix of linear forms in $x_1,\dots,x_4$. The
ideal of $Y$ is generated by the 14 equations
 \[
 M \begin{pmatrix} y_1\\y_2\\y_3\\y_4 \end{pmatrix}=0, \quad
 y_iy_j = (M^\dag)_{ij};
 \]
here $M^\dag$ is the adjoint matrix of signed maximal minors. $Y$ is the
determinantal hypersurface $Y:(\det M=0)\subset\PP^3=\PP(2,2,2,2)$, and
the symmetric matrix $M$ is the resolution of an ample divisorial sheaf
$\Oh_Y(A)$ (singular at the 10 nodes, where $\rank M=2$), such that the
$x_i\in H^0(\Oh_Y(2A)$, $y_i\in H^0(\Oh_Y(3A))$.
 \end{enumerate}
Then $Y$ has a projection of Type~IV from any node $P$ to a codimension~2
complete intersection $X_{6,6}\subset\PP(2,2,2,3,3)$.

Working top down from $Y$, I choose coordinates so that the node
$P=(1,0,0,0)$ in $\PP(2,2,2,2)$, and $y_1,y_2$ are local coordinates at
$P$ (compare \cite{CPR}, 3.4 for local coordinates on quotient singularities
and their fractional divisor of zeros on a blowup). The geometric meaning of
projection is that I blow up $P$ to give an exceptional $-2$-curve $E$, and
consider the blown up surface polarised by the new Weil divisor $A-\half E$.
This eliminates $x_1$ because it does not vanish on $E$; moreover, it
eliminates $y_1,y_2$ because these have weight~3 in the graded ring, but
vanish on $E$ with multiplicity exactly $\half$ (see \cite{CPR}, 3.4).
This is the mark of a Type~IV projection. It is an exercise to make the
appropriate preliminary deduction about $M$ at $P$, and to do the
elimination to get two equations of degree $6$ in the remaining variables.

Now work bottom up from $X$. As in Section~\ref{sec!emb}, one checks that
 \begin{equation}
 \PP^1\to\PP(2,2,2,3,3) \quad\hbox{given by}\quad
 (u,v)\mapsto (u^2,uv,v^2,u^3,v^3)
 \end{equation}
 (or a more general pair of cubics if you prefer) is an embedding, that the
image is contained in the obvious conic $q=x_1x_3-x_2^2$ and in the
3 sextics $f=x_1^3-y_1^2$, $g=x_1x_2x_3-y_1y_2$, $h=x_3^3-y_2^2$, and that
two general sextic linear combinations of $q,f,g,h$ define a quasismooth
K3 surface $X_{6,6}$.

I can resolve $\C[u,v]$ as a module over $\C[\ux,\uy]$ as in (\ref{eq!Gor}),
obtaining
 \begin{equation}
 \C[u,v] \longleftarrow L_0 \xleftarrow{\,P\,} L_1 \xleftarrow{\,M\,} L_2
 \xleftarrow{\,\tP\,} L_3 \ot 0,
 \label{eq!res2}
 \end{equation}
where
 \begin{align*}
 L_0&= A\oplus 2A(-2), \\
 L_1&= 4A(-3)\oplus5A(-4), \\
 L_2&= 5A(-6)\oplus4A(-7), \\
 L_3&= 2A(-9)\oplus A(-10).
 \end{align*}
Here $M$ is a skew $9\times9$ matrix with generic rank~6. Next, following
the argument of (\ref{eq!comps}--\ref{eq!cxs}), I can compose two maps
between complexes to map the Koszul complex of $X$ to
 \begin{equation}
 \renewcommand{\arraystretch}{1.5}
 \begin{matrix}
 A & \ot & A(-6)\oplus A(-6) & \ot & A(-12) & \ot & 0 \\
 \bigcap && \big\downarrow && \big\downarrow \\
 L_0 & \ot & L_1 & \ot & L_2 & \ot & L_3 & \ot & 0,
 \end{matrix}
 \label{eq!cxs2}
 \end{equation}
where the final downarrow defines an element $\mathrm{col}\in L_2(12)$
analogous to (\ref{eq!N24}). The 3 unprojection variables $x_4,y_3,y_4$
in degrees $2,3,3$ correspond to the basis of $L_3$, and as in
(\ref{eq!sp}), the linear relations between them are in the Kustin--Miller
form $\tP\left( \begin{smallmatrix} y_4 \\ y_3 \\ x_4
\end{smallmatrix}\right)=\mathrm{col}$. I believe that the quadratic
relations for $y_3^2,y_3y_4,y_4^2$ can be worked out as in
Section~\ref{sec!quad}.

This example is instructive: although both families of K3 surfaces are
well understood, the equations describing the birational map between them
are very subtle, and would require a lot of computation to elucidate.

 \section{Magma routines}
 \subsection{Embedding $\Ga$ and quasismoothness of $X_4$} \label{ssec!smooth}
 The first thing is to make the polynomial ring $S=\Q[u,v,x,y,z,t]$ and
define the ideal of the graph of the embedding $\PP(2,3)\to\PP(4,5,6,9)$:
 \begin{verbatim}
> Q := Rationals();
> S<u,v,x,y,z,t> := PolynomialRing(Q,[2,3,4,5,6,9]);
> I := Ideal([-x+u^2,-y+u*v,-z+v^2+u^3,-t+u^3*v+v^3]); \end{verbatim}
(you can print $I$ by doing ``I;'' at the prompt, and similarly throughout).

Now eliminate the first two generators to get the ideal of the image $\Ga$:
 \begin{verbatim}
> J := EliminationIdeal(I,2);
> [WeightedDegree(Basis(J)[i]) : i in [1..#Basis(J)]]; \end{verbatim}
answer: $[20,18,19,21,22,23]$. This says that the 1st, 2nd and 3rd basis
elements have weight $20,18,19$. (At this point I used the output of ``J;''
to write (\ref{eq!Ga}).) I make their generic linear combination of weight
24 by hand:
 \begin{verbatim}
> h := Basis(J)[1]; g := Basis(J)[3]; f := Basis(J)[2];
> F0 := x*h + y*g + z*f; IsHomogeneous(F0); \end{verbatim}
answer: true.

 \begin{rmk} \rm
 Taking generic values such as $1,1,1$ for coefficients in a linear system
and jiggling them if necessary is a computer algebra substitute for
Bertini's theorem. Quasismooth is an open condition, so if it ever holds,
we would be infinitely unlucky to happen on a singular guy. In fact
both $F_0=xh+yg-zf$ and $xh+yg$ give singular $X_{24}$.
 \end{rmk}

To work in the polynomial subring $R=k[x,y,z,t]$, the type checking in
Magma insists that I force the elements $x,y,z,t$ into it explicitly 
(otherwise it has no way of knowing that the user wants to identify the
variables in $R$ and $S$ just because they have the same human names).
One way is to set up a homomorphism $\fie\colon S\to R$ taking the
generators $u,v\mapsto0,0$ and $x,y,z,t\mapsto X,Y,Z,T$:
 \begin{verbatim}
> R<X,Y,Z,T> := PolynomialRing(Q,[4,5,6,9]);
> fie := hom< S -> R | 0,0,X,Y,Z,T >;
> F := fie(F0); \end{verbatim}

Now I define the hypersurface $X_{24}:(F=0)\subset\PP(4,5,6,9)$ and ask
the big question:
 \begin{verbatim}
> PP:=Proj(R); X:=Scheme(PP,F);
> IsNonSingular(X); \end{verbatim}
answer: true.

For completeness, I prove that $\PP(2,3)\to\Ga$ is an isomorphism. For
this, I check that monomials of weight $i$ in $\C[X,Y,Z,T]$ map to
polynomials in $u,v$ that span the vector space of polynomials of weight
$i$:
 \begin{verbatim}
> RR<U,V> := PolynomialRing(Q,[2,3]);
> psi := hom< R -> RR | U^2, U*V, U^3+V^2, U^3*V+V^3 >;
> for i in [12..24] do
>   K := ideal< RR | psi(MonomialsOfWeightedDegree(R,i)) >;
>   [m in K : m in MonomialsOfWeightedDegree(RR,i)];
> end for; \end{verbatim}
the last line returns the list of Boolean value (is $\psi(m)\in K$?) for
every monomial in $X,Y,Z,T$ of degree $i\in[12,\dots,24]$. The answer
[true, false, false], [true, true], etc., says that the 2 monomials $U^3$
and $V^2$ in degree~12 fail the test, and every other monomial passes.

 \subsection{Resolution of $\C[u,v]$} \label{ssec!res}
 The problem is to treat $\C[u,v]$ as a module over the ring $\C[x,y,z,t]$,
where the action is as in (\ref{eq!elim}). The following is based on a
Magma routine written by Alan Steel and Gavin Brown. We treat $\C[u,v]$ as
the quotient
 \begin{equation}
 S/(x-u^2,y-uv,z-u^3-v^2,t-u^3v-v^3),
 \end{equation}
where $S=\C[u,v,x,y,z,t]$ (that is, the coordinate ring of the graph of
$i$). It is certainly generated by $1,u,v$ over the big ring $S$ (in fact,
by 1 only), and it is trivial to calculate the matrix of relations:
 \begin{equation}
 \begin{pmatrix}
 u & -1 & 0 \\
 v & 0 & -1 \\
 x-u^2 & 0 & 0 \\
 y-uv & 0 & 0 \\
 z-u^3-v^2 & 0 & 0 \\
 t-u^3v-v^3 & 0 & 0
 \end{pmatrix}
 \begin{pmatrix} 1 \\ u \\ v
 \end{pmatrix} = 0.
 \label{eq!rel1}
 \end{equation}

Define $L\subset F$ to be the submodule of the free module $F=S^3$
generated by the rows of this matrix, and Gr\"obnerise it. In Magma
Gr\"obner is a transitive verb that takes as its object a submodule
of a free module.
 \begin{verbatim}
> Q := Rationals(); // Omit this line at your peril.
> S<u,v,x,y,z,t>:=PolynomialRing(Q,[2,3,4,5,6,9]);
> Free := Module(S,3);
> L:=sub< Free | [ [u, -1, 0], [v, 0, -1], [-u^2+x, 0, 0],
>     [-u*v+y, 0, 0], [-u^3-v^2+z, 0, 0],
>     [-u^3*v-v^3+t, 0, 0] ] >;
> Groebner(L); \end{verbatim}
Gr\"obnerising does many row operations on the matrix in (\ref{eq!rel1}),
to put the dependence on $u$ at the top and on the right, then the same
for $v$. It gives the new basis of $L$:
 \begin{equation}
 \Basis(L) = 
 \begin{pmatrix}
 -y& 0& u \\
 -x& u& 0 \\
 u& -1& 0 \\
 -z& x& v \\
 -y& v& 0 \\
 v& 0& -1 \\
 -xz + y^2& x^2& 0 \\
 x^2& -z& y \\
 -z^2& xz& t \\
 0&xt^2 - y^2z^2& 0 \\
 xt& -yz& 0 \\
 0& -y& x \\
 -z^2t^2& y^2z^3& t^3 \\
 yz& -t& 0 \\
 -t& 0& z
 \end{pmatrix}
 \end{equation}
The rows not involving $u,v$ generate the submodule $N\subset F=S^3$
consisting of all linear relations not involving $u,v$ in the
coefficients.
 \begin{verbatim}
> L1 := [f : f in Basis(L) | &and[Degree(f[i],u) le 0 and
>   Degree(f[i],v) le 0 : i in [1..3] ]];          
> N := sub< Free | L1 >;
> MinimalBasis(N); \end{verbatim}
answer: the minimal basis of $N$ is the matrix (\ref{eq!M}), up to
cosmetic changes.

 \subsection{The K3 Database in Magma} \label{ssec!K3DB}
 Export 2.8 of the computer algebra system Magma \cite{Ma} includes the
first public version of Gavin Brown's K3 surface database, which is
largely a reworking in computer language of the calculations of
Anthony Iano-Fletcher and Selma Alt{\i}nok's Warwick theses \cite{Fl},
\cite{A}. It lists families of polarised K3 surfaces in searchable form.
Further mathematical presentation, and more computer documentation will
be available soon (see \cite{ABR}). Here I only illustrate how to use
the database to find the projections of Type~IV mentioned in the
introduction.
 \begin{verbatim}
> D:=K3Database("t");
> #D; \end{verbatim}
answer: the database $D$ contains 391 elements. The first line defines $D$
to be the K3Database, and requests that Magma refers to the variable in
Hilbert series by its human name $t$, rather than its preferred internal
representation \$.1. (You can put empty brackets $()$ if you like, but if
you omit the brackets, Magma will just sulk and refuse to play.)
 \begin{verbatim}
> #[ X : X in D | Codimension(X) eq 4 ]; \end{verbatim}
answer: 142. Note the logic of Magma's very useful and flexible sequence
construction $\mathtt{[\ ]}$. Here it constructs the sequence of
$X\in D$ such that $\codim X=4$. It runs through all $X\in D$ (the entry
after the colon), evaluates the Boolean value of the statement following
the bar $|$ (is $\codim X=4$?), and returns $X$ (the entry before the
colon) if the answer is yes. We could ask it to return any function of
$X$ by putting $f(X)$ before the colon. Asking for the number of items in
a list $\#[\ ]$ is usually preferable to asking for the list itself, which
may be huge and uninformative. An alternative is to give the list a name,
and ask to see just a bit of it:
 \begin{verbatim}
> cod4 := [ X : X in D | Codimension(X) eq 4 ];
> cod4[1..5]; \end{verbatim}

To play with projections, we first ask Magma to calculate the possible
projection of Type~I, II and~IV of each element of the database $D$ from
each of its singularities:
 \begin{verbatim}
> Centres(~D);
> Cod4Type4 := [ X : X in D | Codimension(X) eq 4 \
>  and &or[ p[2] eq [4] : p in Centres(X)] ];\end{verbatim}
The last line here defines Cod4Type4 to be the sequence of $X\in D$
having $\codim X=4$ and such that at least one of the centres $p$ of
$X$ has Type~IV. The internal $[\ ]$ constructs a sequence of Boolean
values (for $p$ a singularity of $X$, is the projection from $p$ of
Type~IV?) and the $\&\mathrm{or}[\ ]$ evaluates ``or''. To write this code
you need to know a little about how the database stores its data. You can
find out by interrogating it. For example,
 \begin{verbatim}
> X := D[293]; Centres(X); \end{verbatim}
tells you that No.~293 in $D$ has a singularity of type $\frac15(2,3)$ with
a projection of Type~I, etc., and you figure out that each centre $p$ is
stored as a sequence of entries of which the second $p[2]$ is the type of
the projection. I ask how many elements does Cod4Type4 contain, and what
are the first two elements?
 \begin{verbatim}
> #Cod4Type4;
> Cod4Type4[1..2]; \end{verbatim}
answer: 5 elements. For the first, see below. The second is the example
treated in Sections~\ref{sec!emb}--\ref{sec!quad} above. All five of these
K3s are from Selma Alt{\i}nok's thesis. For historical reasons, the
database lists the codimension~4 K3s as Altinok4$(-24)$ through to
Altinok4(121), missing Nos.~0, 76, 81, 107. The 5 elements of Cod4Type4
are listed as Altinok4(84), (111), (53), $(-11)$ and $(-6)$.

 \begin{exa} \rm
 The first element of the above list is Altinok4(84), the codimension~4
K3 surface $Y\subset\PP(5,5,6,7,8,9,11)$ with basket
$2\times\frac15(2,3),\frac1{11}(5,6)$ and Hilbert numerator
 \[
 1-t^{14}-t^{15}-2t^{12}-t^{17}-2t^{18}-t^{19}-t^{20}+\cdots
 \]
The database says that this has a projection of Type~IV to a surface
having the same numerical type as $X_{66}\subset\PP(5,6,22,33)$, but
containing a copy of $\PP(2,3)$. This could happen because $X_{66}$ is
actually a hypersurface of $\PP(5,6,22,33)$ containing the image of a
general embedding $\PP(2,3)\into\PP(5,6,22,33)$; or $X_{66}$ could
itself be a bit degenerate, for example, with the curve $\PP(2,3)$ as a
base locus of some monogonal linear system. It would be interesting to
see this and the other cases in the list worked out in the spirit of
this paper.
 \end{exa}

 \begin{rmk} \rm
 We expect subsequent exports of Magma to contain substantial improvements
to the internal structures of the K3 database, its contents, and its user
interface. Whereas at present we use Magma mainly as a book-keeping device
to handle large quantities of combinatorial data, we hope that in future
we can use it to automate many of the computations in commutative and
homological algebra involved in the geometry of K3 surfaces and Fano
3-folds, including the calculations above and those in \cite{CPR}
and \cite{Ki}.
 \end{rmk}

\bigskip\noindent
Miles Reid, \\
Math Inst., Univ.\ of Warwick, \\
Coventry CV4 7AL, England \\
e-mail: miles@maths.warwick.ac.uk \\
web: www.maths.warwick.ac.uk/$\!\scriptstyle\sim$miles

\end{document}